\ifx\documentclass\undefined
\documentstyle[11pt]{article}
\else
\documentclass[11pt]{article}
\usepackage{latexsym}
\usepackage{amsfonts}
\usepackage{amsmath}
\usepackage{amssymb}
\fi

\author{ K\'aroly Bezdek 
\thanks{Partially supported by a Natural Sciences and 
Engineering Research Council of Canada Discovery Grant}}

\font\tenBbb=msbm10 at 12pt         \font\sevenBbb=msbm9    \font\fiveBbb=msbm7
\newfam\Bbbfam
\textfont\Bbbfam=\tenBbb \scriptfont\Bbbfam=\sevenBbb
\scriptscriptfont\Bbbfam=\fiveBbb

\def\kkk{\null\hfill $\Box$\smallskip}

\newcommand{\proof}{{\noindent\bf Proof:{\ \ }}}

\newtheorem{theorem}{Theorem}[section]
\newtheorem{lemma}[theorem]{Lemma}

\newtheorem{cor}[theorem]{Corollary}

\title{On a strong version of the Kepler conjecture 
}

\begin{document}

\maketitle

\date

\begin{abstract}
We raise and investigate the following problem that one can regard as a very close relative of the densest sphere packing problem. If the Euclidean $3$-space is partitioned into convex cells each containing a unit ball, how should the shapes of the cells be designed to minimize the average surface area of the cells?  In particular, we prove that the average surface area in question is always at least $\frac{24}{\sqrt{3}}=13.8564...$. 
\end{abstract}

\medskip

\section{Introduction}

The central problem that we raise in this paper can be phrased informally as follows: {\it if the Euclidean $3$-space is partitioned into convex cells each containing a unit ball, how should the shapes of the cells be designed to minimize the average surface area of the cells}?
In order to state our problem in more precise terms we proceed as follows. Let $\cal T$ be a tiling of the $3$-dimensional Euclidean space $\mathbb{E}^{3}$
into convex polyhedra $\mathbf{P}_i, i=1, 2, \dots$ each containing a unit ball say, $\mathbf{P}_i$ containing the closed $3$-dimensional ball $\mathbf{B}_i$ centered at 
the point $\mathbf{o}_i$ having radius $1$ for $i=1, 2, \dots$. Also, we assume that there is a finite upper bound for the diameters of the convex cells in $\cal T$, i.e., $\sup\{{\rm diam}(\mathbf{P}_i)|i=1, 2, \dots\}<\infty $, where ${\rm diam}(\cdot )$ denotes the diameter of the corresponding set. In short, we say that $\cal T$ is a {\it normal tiling} of  $\mathbb{E}^{3}$ with the underlying packing $\cal P$ of the unit balls $\mathbf{B}_i, i=1, 2, \dots$. Then we define the (lower) {\it average surface area} $\underline{s}({\cal T})$ of the cells in $\cal T$ as follows:
$$
\underline{s}({\cal T}):=\liminf_{L\to\infty}\frac{\sum_{\{i|\mathbf{B}_i\subset \mathbf{C}_{L}\}}^{}{\rm sarea}(\mathbf{P}_i\cap \mathbf{C}_{L})}{{\rm card}\{i|\mathbf{B}_i\subset \mathbf{C}_{L}\}},
$$
where $ \mathbf{C}_{L}$ denotes the cube centered at the origin $\mathbf{o}$ with edges parallel to the coordinate axes of $\mathbb{E}^{3}$ and having edge length ${L}$ furthermore, ${\rm sarea}(\cdot )$ and ${\rm card}(\cdot )$ denote the surface area and cardinality of the corresponding sets.  (We note that it is rather straightforward to show that  $\underline{s}({\cal T})$ is independent from the choice of the coordinate system of $\mathbb{E}^{3}$.)

There is very natural way to generate a large family of normal tilings. Namely, let ${\cal P}_R$ be an arbitrary packing of unit balls in $\mathbb{E}^{3}$ with the property
that each closed ball of radius $R$ in $\mathbb{E}^{3}$ contains the center of at least one unit ball in ${\cal P}_R$. Recall that the Voronoi cell of a unit ball in ${\cal P}_R$ is the set of points that are not farther away from the center of the given ball than from any other ball's center.  It is well known that the Voronoi cells in question form a tiling of $\mathbb{E}^{3}$ (for more details see \cite{Ro64}). Furthermore, the Voronoi tiling obtained in this way, is going to be a normal one because each Voronoi cell is contained in the closed ball of radius $R$ concentric to the unit ball of the given Voronoi cell and therefore the diameter of each Voronoi cell is at most $2R$. Also, we recall here the {\it strong dodecahedral conjecture} of \cite{Bez00}: the surface area of every (bounded) Voronoi cell in a packing of unit balls  is at least that of a regular dodecahedron of inradius $1$, i.e., it is at least $16.6508\dots$. After a sequence of partial results obtained in \cite{Bez00}, \cite{BeDa05}, and \cite{AmFo06} (proving the lower bounds $16.1433\dots$, $16.1445\dots$, and $16.1977\dots$), just very recently, Hales \cite{Ha11} has announced a computer assisted proof of the strong dodecahedral conjecture.

In the second half of this paper, by adjusting Kert\'esz's volume estimation technique (\cite{Ke88}) to our problem on estimating surface area and making the necessary modifications, we give a proof of the following inequality. 

\begin{theorem}\label{Bezdek-estimate}
Let $\cal T$ be an arbitrary normal tiling of $\mathbb{E}^{3}$. Then the average surface area
of the cells in $\cal T$ is always at least $\frac{24}{\sqrt{3}}$, i.e.,
$$
\underline{s}({\cal T})\ge \frac{24}{\sqrt{3}}=13.8564...\ .
$$
\end{theorem} 

Most likely the lower bound in Theorem~\ref{Bezdek-estimate} can be improved further however, any such improvement would require additional new ideas. In particular, recall that in the face-centered cubic lattice packing of unit balls in $\mathbb{E}^{3}$, when each ball is touched by $12$ others, the Voronoi cells of the unit balls are regular rhombic dodecahedra of inradius $1$ and of surface area $12\sqrt{2}$ (for more details on the geometry involved see \cite{FeTo64}). Thus, it is immediate to raise the following question: {\it prove or disprove that if $\cal T$ is an arbitrary normal tiling of  $\mathbb{E}^{3}$, then}

\begin{equation}\label{Bezdek-problem} 
\underline{s}({\cal T})\ge 12\sqrt{2}=16.9705... \ . 
\end{equation}

Let us mention that an affirmative answer to (\ref{Bezdek-problem}) for the family of Voronoi tilings of unit ball packings would imply the Kepler conjecture. As is well known, the Kepler conjecture has been proved by Hales in a sequence of celebrated papers (\cite{Ha05}, \cite{HaFe06}, \cite{Ha06-III},  \cite{Ha06-IV}, and  \cite{Ha06-VI}) concluding that the density of any unit ball packing in $\mathbb{E}^{3}$ is at most $\frac{\pi}{\sqrt{18}}$. Indeed, if $\underline{s}({\cal T})\ge 12\sqrt{2}$ were true for the Voronoi tilings $\cal T$ of unit ball packings $\cal P$ in $\mathbb{E}^{3}$, then based on the obvious inequalities
$$\sum_{\{i|\mathbf{B}_i\subset \mathbf{C}_{L}\}}^{}{\rm vol}(\mathbf{P}_i\cap \mathbf{C}_{L})\le{\rm vol}(\mathbf{C}_{L})\ \ {\rm and}\ \ 
\frac{1}{3}{\rm sarea}(\mathbf{P}_i\cap \mathbf{C}_{L})\le {\rm vol}(\mathbf{P}_i\cap \mathbf{C}_{L}),$$
(where ${\rm vol}(\cdot )$ denotes the volume of the corresponding set) we would get that the (upper) {\it density} 
$
\overline{\delta}({\cal P}):=\limsup_{L\to\infty}\frac{\frac{4\pi}{3}{\rm card}\{i|\mathbf{B}_i\subset \mathbf{C}_{L}\}}{{\rm vol}(\mathbf{C}_{L})} 
$ 
of the packing ${\cal P}$ must satisfy the inequality 
$$\overline{\delta}({\cal P})\le\limsup_{L\to\infty}\frac{\frac{4\pi}{3}{\rm card}\{i|\mathbf{B}_i\subset \mathbf{C}_{L}\}}{ \sum_{\{i|\mathbf{B}_i\subset \mathbf{C}_{L}\}}^{}{\rm vol}(\mathbf{P}_i\cap \mathbf{C}_{L})  }
$$
$$
\le\limsup_{L\to\infty}\frac{4\pi{\rm card}\{i|\mathbf{B}_i\subset \mathbf{C}_{L}\}}{ \sum_{\{i|\mathbf{B}_i\subset \mathbf{C}_{L}\}}^{}{\rm sarea}(\mathbf{P}_i\cap \mathbf{C}_{L})  }
 =\frac{4\pi}{\underline{s}({\cal T})}\le \frac{\pi}{\sqrt{18}}
$$ 
Thus, one could regard the affirmative version of (\ref{Bezdek-problem}), stated for the Voronoi tilings of unit ball packings, as a {\it strong version of the Kepler conjecture}.

As an additional observation we mention that an affirmative answer to (\ref{Bezdek-problem}) would imply also the rhombic dodecahedral conjecture of \cite{Be00}. According to that conjecture the surface area of any $3$-dimensional parallelohedron of inradius at least $1$ (i.e., the surface area of any convex polyhedron containing a unit ball and having a family of translates tiling $\mathbb{E}^{3}$) is at least as large as $12\sqrt{2}=16.9705...$.

  Last but not least, it is very tempting to further relax the conditions in our original problem by replacing convex cells with cells that are measurable and have measurable boundaries and ask the following more general question: {\it if the Euclidean $3$-space is partitioned into cells each containing a unit ball, how should the shapes of the cells be designed to minimize the average surface area of the cells}? One can regard this question as a {\it foam problem}, in particular, as a relative of Kelvin's foam problem (on partitioning $\mathbb{E}^{3}$ into unit volume cells with minimum average surface area) since foams are simply tilings of space that try to minimize surface area. Although foams are well studied (see
the relevant sections of the highly elegant book \cite{Mo09} of Morgan), it is far not clear what would be a good candidate for the proper minimizer in the foam question just raised. As a last note we mention that Brakke \cite{Br11}, by properly modifying the Williams foam, has just obtained a partition of the Euclidean $3$-space into congruent cells each containing a unit ball and having surface area $16.95753<12\sqrt{2}=16.9705...$.

\section{Proof of Theorem~\ref{Bezdek-estimate}}

First, we prove the following ``compact'' version of Theorem~\ref{Bezdek-estimate}. It is  also a surface area analogue of the volume estimating theorem in  \cite{Ke88}. 

\begin{theorem}\label{compact-Bezdek-estimate}
If the cube $ \mathbf{C}$ is partitioned into the convex cells $ \mathbf{Q}_1,  \mathbf{Q}_2, \dots ,$ $ \mathbf{Q}_n$ each containing a unit ball in $\mathbb{E}^{3}$, then the sum of the surface areas of the $n$ convex cells is at least $\frac{24}{\sqrt{3}}n$, i.e.,
$$
\sum_{i=1}^n{\rm sarea}(\mathbf{Q}_i)\ge \frac{24}{\sqrt{3}}n\ .
$$
\end{theorem}

\proof
Let $E(\mathbf{Q}_i)$ denote the family of the edges of the convex polyhedron $\mathbf{Q}_i$ and let ${\rm ecurv}(\mathbf{Q}_i):=\sum_{e\in E(\mathbf{Q}_i)}L(e)\tan\frac{\alpha_e}{2}$ be the so-called {\it  edge curvature} of $ \mathbf{Q}_i$, where $L(e)$ denotes the length of the edge $e\in E(\mathbf{Q}_i)$ and $\alpha_e$ is the angle between the outer normal vectors of the two faces of   $ \mathbf{Q}_i$ meeting along the edge $e$, $1\le i\le n$. It is well known that the Brunn-Minkowski inequality implies the following inequality (for more details we refer the interested reader to p. 287 in \cite{FeTo64}):
\begin{equation}\label{inequality-1}
{\rm sarea}^2(\mathbf{Q}_i)\ge 3 {\rm vol}(\mathbf{Q}_i){\rm ecurv}(\mathbf{Q}_i)\ .
\end{equation}
Also, it will be more proper for us to use the inner dihedral angles $\beta_e:=\pi-\alpha_e$ and the relevant formula
\begin{equation}\label{equality-1}
{\rm ecurv}(\mathbf{Q}_i)=\sum_{e\in E(\mathbf{Q}_i)}L(e)\cot\frac{\beta_e}{2}\ .
\end{equation}
As, by assumption, $ \mathbf{Q}_i$ contains a unit ball therefore
\begin{equation}\label{inequality-2}
{\rm vol}(\mathbf{Q}_i)\ge \frac{1}{3}{\rm sarea}(\mathbf{Q}_i)\ .
\end{equation}
Hence, (\ref{inequality-1}), (\ref{equality-1}), and (\ref{inequality-2}) imply in a straightforward way that
\begin{equation}\label{inequality-3}
{\rm sarea}(\mathbf{Q}_i)\ge\sum_{e\in E(\mathbf{Q}_i)}L(e)\cot\frac{\beta_e}{2}
\end{equation}
holds for all $1\le i\le n$.

Now, let $s\subset\mathbf{C}$ be a closed line segment along which exactly $k$ members of the family $\{ \mathbf{Q}_1,  \mathbf{Q}_2, \dots ,$ $ \mathbf{Q}_n\}$ meet having inner dihedral angles $\beta_1, \beta_2, \dots , \beta_k$. There are the following three possibilities:

{\it (a)}  $s$ is on an edge of the cube $ \mathbf{C}$;

{\it (b)}  $s$ is in the relative interior either of a face of  $ \mathbf{C}$ or of a face of a convex cell in the family $\{ \mathbf{Q}_1,  \mathbf{Q}_2, \dots ,$ $ \mathbf{Q}_n\}$;

{\it (c)}  $s$ is in the interior of $ \mathbf{C}$ and not in the relative interior of any face of any convex cell in the family $\{ \mathbf{Q}_1,  \mathbf{Q}_2, \dots ,$ $ \mathbf{Q}_n\}$.

In each of the above cases we can make the following easy observations:

{\it (a)}  $\beta_1+ \beta_2+ \dots + \beta_k=\frac{\pi}{2}$ with $k\ge 1$;

{\it (b)}  $\beta_1+ \beta_2+ \dots + \beta_k=\pi$ with $k\ge 2$;

{\it (c)}  $\beta_1+ \beta_2+ \dots + \beta_k=2\pi$ with $k\ge 3$.

As $y=\cot x$ is convex and decreasing over the interval $0<x\le\frac{\pi}{2}$ therefore the following inequalities must hold:

{\it (a)} $\cot\frac{\beta_1}{2}+\cot\frac{\beta_2}{2}+\dots +\cot\frac{\beta_k}{2}\ge k\cot\frac{\pi}{4k}\ge k$;

{\it (b)} $\cot\frac{\beta_1}{2}+\cot\frac{\beta_2}{2}+\dots +\cot\frac{\beta_k}{2}\ge k\cot\frac{\pi}{2k}\ge k$;

{\it (c)} $\cot\frac{\beta_1}{2}+\cot\frac{\beta_2}{2}+\dots +\cot\frac{\beta_k}{2}\ge k\cot\frac{\pi}{k}\ge \frac{1}{\sqrt{3}}k$.

In short, the following inequality holds in all three cases:
\begin{equation}\label{inequality-4}
\cot\frac{\beta_1}{2}+\cot\frac{\beta_2}{2}+\dots +\cot\frac{\beta_k}{2}\ge \frac{1}{\sqrt{3}}k\ .
\end{equation}

Thus, by adding together the inequalities (\ref{inequality-3}) for all $1\le i\le n$ and using (\ref{inequality-4}) we get that

\begin{equation}\label{inequality-5}
\sum_{i=1}^n{\rm sarea}(\mathbf{Q}_i)\ge\frac{1}{\sqrt{3}}\sum_{i=1}^{n}\sum_{e\in E(\mathbf{Q}_i)}L(e)\ .
\end{equation}

Finally, recall the elegant theorem of Besicovitch and Eggleston \cite{BeEg57} claiming that the total edge length  of any convex polyhedron containing a unit ball in $\mathbb{E}^{3}$ is always at least as large as the total edge length of a cube circumscribed a unit ball. This implies that
\begin{equation}\label{inequality-6}
\sum_{e\in E(\mathbf{Q}_i)}L(e)\ge 24
\end{equation}
holds for all $1\le i\le n$. Hence, (\ref{inequality-5}) and (\ref{inequality-6}) finish the proof of Theorem~\ref{compact-Bezdek-estimate}.      \kkk

Second, we take a closer look of the given normal tiling $\cal T$ defined in details in the first Section of this paper and using Theorem~\ref{compact-Bezdek-estimate} we give a proof of Theorem~\ref{Bezdek-estimate}.

By assumption $D:=\sup\{{\rm diam}(\mathbf{P}_i)|i=1, 2, \dots\}<\infty$. Thus,  clearly each closed ball of radius $D$ in $\mathbb{E}^{3}$ contains at least one of the convex polyhedra $\mathbf{P}_i, i=1, 2, \dots$ (forming the tiling $\cal T$ of $\mathbb{E}^{3}$). Now, let $ \mathbf{C}_{L_N}, N=1, 2, \dots$ be an arbitrary sequence of cubes centered at the origin $\mathbf{o}$ with edges parallel to the coordinate axes of $\mathbb{E}^{3}$ and having edge length ${L_N}, N=1, 2, \dots$ with $\lim_{N\to\infty}L_n=\infty$. It follows that
\begin{equation}\label{inequality-7}
0<\liminf_{N\to\infty}\frac{\frac{4\pi}{3}{\rm card}\{i|\mathbf{B}_i\subset \mathbf{C}_{L_N}\}}{{\rm vol}(\mathbf{C}_{L_N})}\le
\limsup_{N\to\infty}\frac{\frac{4\pi}{3}{\rm card}\{i|\mathbf{B}_i\subset \mathbf{C}_{L_N}\}}{{\rm vol}(\mathbf{C}_{L_N})}<1 . 
\end{equation}
Note that clearly
\begin{equation}\label{inequality-8}
\frac{{\rm card}\{i|\mathbf{P}_i\cap{\rm bd}\mathbf{C}_{L_N}\neq\emptyset \}}{{\rm card}\{i|\mathbf{B}_i\subset \mathbf{C}_{L_N}\}}
\le
\frac{\big({\rm vol}(\mathbf{C}_{L_N+2D})- {\rm vol}(\mathbf{C}_{L_N-2D})\big)   {\rm vol}(\mathbf{C}_{L_N})}{{\rm vol}(\mathbf{C}_{L_N})  \frac{4\pi}{3} {\rm card}\{i|\mathbf{B}_i\subset \mathbf{C}_{L_N}\} }
\end{equation}
moreover,
\begin{equation}\label{inequality-9}
\lim_{N\to\infty}\frac{{\rm vol}(\mathbf{C}_{L_N+2D}) - {\rm vol}(\mathbf{C}_{L_N-2D})  }{{\rm vol}(\mathbf{C}_{L_N})}=0.
\end{equation}
Thus, (\ref{inequality-7}), (\ref{inequality-8}), and (\ref{inequality-9}) imply in a straightforward way that
\begin{equation}\label{inequality-10}
\lim_{N\to\infty}
\frac{{\rm card}\{i|\mathbf{P}_i\cap{\rm bd}\mathbf{C}_{L_N}\neq\emptyset \}}{{\rm card}\{i|\mathbf{B}_i\subset \mathbf{C}_{L_N}\}}=0 \ .
\end{equation}
Moreover, (\ref{inequality-3}) yields that
\begin{equation}\label{inequality-11}
{\rm sarea}(\mathbf{P}_i)\ge{\rm ecurv}(\mathbf{P}_i)=\sum_{e\in E(\mathbf{P}_i)}L(e)\cot\frac{\beta_e}{2}
\end{equation}
holds for all $i=1, 2, \dots$.
As a next step, using
\begin{equation}\label{inequality-12}
{\rm sarea}(\mathbf{P}_i)={\rm sarea}\left( {\rm bd} (\mathbf{P}_i\cap\mathbf{C}_{L})\setminus{\rm bd}\mathbf{C}_{L}\right)+\delta_i
\end{equation}
and
\begin{equation}\label{inequality-13}
{\rm ecurv}(\mathbf{P}_i)\ge \sum_{e\in E\left({\rm bd}(\mathbf{P}_i\cap\mathbf{C}_{L})\setminus{\rm bd}\mathbf{C}_{L}\right)}L(e)\cot\frac{\beta_e}{2}
\end{equation}
(with ${\rm bd}(\cdot )$ denoting the boundary of the corresponding set) we obtain the following from (\ref{inequality-11}):
\begin{equation}\label{inequality-14}
{\rm sarea}\left( {\rm bd}(\mathbf{P}_i\cap\mathbf{C}_{L})\setminus{\rm bd}\mathbf{C}_{L}\right)+\delta_i
\ge
\sum_{e\in E\left({\rm bd}(\mathbf{P}_i\cap\mathbf{C}_{L})\setminus{\rm bd}\mathbf{C}_{L}\right)}L(e)\cot\frac{\beta_e}{2}\ ,
\end{equation}
where clearly $0\le\delta_i\le{\rm sarea}(\mathbf{P}_i)$.
Hence, (\ref{inequality-14}) combined with (\ref{inequality-4}) yields 
\begin{cor}\label{corollary-of-proof}
$$
f(L):=\sum_{\{i|{\rm int}\mathbf{P}_i\cap\mathbf{C}_{L}\neq\emptyset\}}{\rm sarea}\big( {\rm bd}(\mathbf{P}_i\cap\mathbf{C}_{L})\setminus{\rm bd}\mathbf{C}_{L}\big)+ \sum_{\{i|\mathbf{P}_i\cap{\rm bd}\mathbf{C}_{L}\neq\emptyset\}}\delta_i 
$$
$$
\ge g(L):=\frac{1}{\sqrt{3}}\sum_{\{i|{\rm int}\mathbf{P}_i\cap\mathbf{C}_{L}\neq\emptyset\}}\bigg(\sum_{e\in E\left({\rm bd}(\mathbf{P}_i\cap\mathbf{C}_{L})\setminus{\rm bd}\mathbf{C}_{L}\right)}L(e)\bigg)
$$
\end{cor}
Now, it is easy to see that 
\begin{equation}\label{inequality-15}
f(L)=\Delta (L)+\sum_{\{i|\mathbf{B}_i\subset \mathbf{C}_{L}\}}^{}{\rm sarea}(\mathbf{P}_i\cap \mathbf{C}_{L})\ ,
\end{equation}
where $0\le\Delta(L)\le 2\sum_{\{i|\mathbf{P}_i\cap{\rm bd}\mathbf{C}_{L}\neq\emptyset\}}{\rm sarea}(\mathbf{P}_i)$.

Moreover, (\ref{inequality-6}) implies that
\begin{equation}\label{inequality-16}
g(L)\ge -\overline{\Delta}(L)+\sum_{\{i|\mathbf{B}_i\subset \mathbf{C}_{L}\}}\frac{24}{\sqrt{3}}\ ,
\end{equation}
where $0\le\overline{\Delta}(L)\le\sum_{\{i|\mathbf{P}_i\cap{\rm bd}\mathbf{C}_{L}\neq\emptyset\}}\sum_{e\in E(\mathbf{P}_i)}L(e) $.

\begin{lemma}\label{sup-area-edge-length}
$$A:=\sup\{{\rm sarea}(\mathbf{P}_i)|i=1, 2, \dots\}<\infty$$ and $$E:=\sup\{\sum_{e\in E(\mathbf{P}_i)}L(e)|i=1, 2, \dots\}<\infty \ .$$
\end{lemma}

\proof
As $D=\sup\{{\rm diam}(\mathbf{P}_i)|i=1, 2, \dots\}<\infty$ therefore according to Jung's theorem (\cite{De85}) each $\mathbf{P}_i $ is contained in a closed ball of radius $\sqrt{\frac{3}{8}}D$  in $\mathbb{E}^{3}$. Thus, $A\le\frac{3}{2}\pi D^2<\infty$.

For a proof of the other claim recall that $\mathbf{P}_i$ contains the unit ball $\mathbf{B}_i$ centered at $\mathbf{o}_i$. If the number of faces of $\mathbf{P}_i$ is $f_i$, then $\mathbf{P}_i$ must have at least $f_i$ neighbours (i.e., cells of $\cal T$ that have at least one point in common with $\mathbf{P}_i$) and as each neighbour is contained in the closed $3$-dimensional ball of radius $2D$ centered at $\mathbf{o}_i$ therefore the number of neighbours of $\mathbf{P}_i$ is at most $(2D)^3-1$ and so, $f_i\le 8D^3-1$. (Here, we have used the fact that each neighbour contains a unit ball and therefore its volume is larger than $\frac{4\pi}{3}$.) Finally, Euler's formula implies that the number of edges of $\mathbf{P}_i$ is at most $3f_i-6\le 24D^3-9$. Thus, $ E\le 24D^4-9D<\infty$ (because the length of any edge of $\mathbf{P}_i$ is at most $D$).  \kkk

Thus, Corollary~\ref{corollary-of-proof}, (\ref{inequality-15}), (\ref{inequality-16}), and Lemma~\ref{sup-area-edge-length} imply the following inequality in a straightforward way.

\begin{cor}\label{main-corollary}
$$
\frac{2A {\rm card}{\{i|\mathbf{P}_i\cap{\rm bd}\mathbf{C}_{L}\neq\emptyset\}} +\sum_{\{i|\mathbf{B}_i\subset \mathbf{C}_{L}\}}^{}{\rm sarea}(\mathbf{P}_i\cap \mathbf{C}_{L})}{ {\rm card}\{i|\mathbf{B}_i\subset \mathbf{C}_{L}\} }
$$
$$
\ge\frac{ -E {\rm card}{\{i|\mathbf{P}_i\cap{\rm bd}\mathbf{C}_{L}\neq\emptyset\}} + \sum_{\{i|\mathbf{B}_i\subset \mathbf{C}_{L}\}}\frac{24}{\sqrt{3}} }{ {\rm card}\{i|\mathbf{B}_i\subset \mathbf{C}_{L}\}  }\ .
$$
\end{cor}
Finally, Corollary~\ref{main-corollary} and (\ref{inequality-10}) yield that
\begin{equation}\label{inequality-17}
\liminf_{N\to\infty}\frac{\sum_{\{i|\mathbf{B}_i\subset \mathbf{C}_{L_N}\}}^{}{\rm sarea}(\mathbf{P}_i\cap \mathbf{C}_{L_N})}{{\rm card}\{i|\mathbf{B}_i\subset \mathbf{C}_{L_N}\}}
\ge\frac{ 24 }{\sqrt{3}  }\ ,
\end{equation}
finishing the proof of Theorem~\ref{Bezdek-estimate}.

\vspace{1cm}

\medskip

\noindent
K\'aroly Bezdek
\newline
Department of Mathematics and Statistics, University of Calgary, Canada,
\newline
Department of Mathematics, University of Pannonia, Veszpr\'em, Hungary,
\newline
and
\newline
Institute of Mathematics, E\"otv\"os University, Budapest, Hungary.
\newline
{\sf E-mail: bezdek@math.ucalgary.ca}


\begin{thebibliography}{GGM}

\bibitem{AmFo06}G. Ambrus and F. Fodor,  A new lower bound on the surface area of a Voronoi polyhedron, \emph{Period. Math. Hungar.} \textbf{53/1-2} (2006), 45--58.

\bibitem{BeEg57}A. S. Besicovitch and H. G. Eggleston, The total length of the edges of a polyhedron, \emph{Quart. J. Math. Oxford Ser.} \textbf{2/ 8} (1957), 172--190.

\bibitem{Bez00}K. Bezdek, On a stronger form of Rogers's lemma and the minimum surface area of Voronoi cells in unit ball packings, \emph{J. Reine Angew. Math.} \textbf{518} (2000), 131--143.

\bibitem{Be00}K. Bezdek, On rhombic dodecahedra, \emph{Contributions
to Alg. and Geom.} \textbf{41/2} (2000), 411--416.

\bibitem{BeDa05}K. Bezdek and E. Dar\'oczy-Kiss, Finding the best face on a Voronoi polyhedron
Ñ the strong dodecahedral conjecture revisited, \emph{Monatsh. Math.} \textbf{145/3} (2005), 191--206.

\bibitem{Br11}K. Brakke, personal communication (Nov., 2011).

\bibitem{De85}B. V. Dekster,  An extension of Jung's theorem, \emph{Israel J. Math.} \textbf{50/3}
(1985), 169--180. 

\bibitem{FeTo64}L. Fejes T\'oth, Regular Figures, \emph{Pergamon Press}, New York, 1964.

\bibitem{Ha05}T. C. Hales, A proof of the Kepler conjecture, \emph{Ann. of Math.}, \textbf{162/3} (2005), 1065--1185.

\bibitem{HaFe06}T. C. Hales and S. P. Ferguson, A formulation of the Kepler conjecture, \emph{Discrete Comput. Geom.} \textbf{36/1} (2006), 21--69.

\bibitem{Ha06-III}T. C. Hales, Sphere packings III, Extremal cases, \emph{Discrete Comput. Geom.} \textbf{36/1} (2006), 71--110.

\bibitem{Ha06-IV}T. C. Hales, Sphere packings IV, Detailed bounds, \emph{Discrete Comput. Geom.} \textbf{36/1} (2006), 111--166.

\bibitem{Ha06-VI}T. C. Hales, Sphere packings VI, Tame graphs and linear programs, \emph{Discrete Comput. Geom.} \textbf{36/1} (2006), 205--265. 

\bibitem{Ha11}T. C. Hales, The strong dodecahedral conjecture and Fejes T\'oth's contact conjecture, \emph{arXiv}:1110.0402v1 [math.MG] (2011), 1--11.

\bibitem{Ke88} G. Kert\'esz, On totally separable packings of equal balls, \emph{Acta Math. Hungar.} \textbf{51/3-4} (1988), 363--364.

\bibitem{Mo09}F. Morgan, Geometric Measure Theory - A Beginner's Guide, Fourth edition, \emph{Elsevier/Academic Press}, Amsterdam, 2009.

\bibitem{Ro64}C. A. Rogers, Packing and Covering, \emph{Cambridge University Press}, Cambridge, 1964.



\end{thebibliography}
\end{document}